\documentclass[11pt]{article}

\usepackage{paralist}
\usepackage[margin=3.0cm]{geometry} 
\usepackage{graphicx} 

\usepackage{amssymb, amsfonts, amsmath, amscd, amsthm} 
\usepackage{mathtools} 
\usepackage{titlesec} 
\usepackage{csquotes} 
\usepackage[shortlabels]{enumitem} 
\usepackage{tikz-cd} 
\usepackage{wrapfig} 
\usepackage{epigraph} 
\usepackage{hyperref} 
\usepackage{todonotes} 
\usepackage{comment}
\usepackage{mathrsfs} 

\hypersetup{
     colorlinks   = true,
     citecolor    = green
}

\numberwithin{equation}{section}

\titleformat*{\section}{\Large \scshape\center}
\titleformat*{\subsection}{\fontsize{14}{14} \sffamily}

\theoremstyle{plain}
\newtheorem{theorem}{Theorem}[section]
\newtheorem*{theorem*}{Theorem}

\newtheorem{lemma}[theorem]{Lemma}
\newtheorem{proposition}[theorem]{Proposition}

\theoremstyle{definition}
\newtheorem{definition}[theorem]{Definition}
\newtheorem*{definition*}{Definition}
\newtheorem{example}[theorem]{Example}

\theoremstyle{remark}
\newtheorem*{remark}{Remark}

\begin{document}
\pagenumbering{gobble}
\title{\huge{A Brief Introduction to \\ the Feichtinger Algebra $\mathbf{
S}_{0}(\mathbb{R})$}}
\author{Eirik Berge}
\date{}
\maketitle
\pagenumbering{arabic}

\begin{center}
    \textit{Dedicated to Hans Georg Feichtinger on the occasion of his 70th birthday.}
\end{center}

\section{Introduction}
In recent decades, time-frequency analysis has gone through a transformation. Having its roots in engineering disciplines, it may come as a surprise that time-frequency analysis is of increasing interest for mathematicians. When asked to give a seminar talk on time-frequency analysis for mathematicians, I desperately needed to answer the question:
\begin{center}
\textit{What single concept in time-frequency analysis \\ should every mathematician know about?}
\end{center} 
Topics such as Gabor frames, the Wigner distribution, and Weyl quantization are intriguing, but are primarily motivated by applied phenomena. Gradually, I realized that the \textit{Feichtinger algebra} $\mathbf{S}_{0}$ works as a gateway to almost all other aspects of time-frequency analysis. On top of this, $\mathbf{S}_{0}$ is a intriguing and beautiful piece of modern mathematics. \par 
This note is based on the seminar talk I gave and highlights the key aspects of the Feichtinger algebra $\mathbf{S}_{0}(\mathbb{R})$ on the real line. There already exists an excellent comprehensive survey of the Feichtinger algebra \cite{mads18} aimed at an audience comfortable with locally compact abelian groups. However, I wrote this note as the literature is lacking a brief pedestrian introduction to the Feichtinger algebra for non-experts.

\subsubsection*{Brief Historical Account}
The Feichtinger algebra $\mathbf{S}_{0}$ was first presented in 1979 at the workshop \textit{Internationale Arbeitstagung über Topologische Gruppen und Gruppenalgebren} by Hans Georg Feichtinger. It was subsequently developed in the papers \cite{feichtinger_original} and \cite{losert80} where the defining feature was a minimality requirement. The Feichtinger algebra is a special case of a \textit{modulation space}, witch originally was introduced in the technical report \cite{modulation_space_report}. Later, Feichtinger and Gr\"{o}chenig developed \textit{coorbit spaces} in \cite{feichtinger1988unified}, where the modulation spaces are prominent examples. The popularity of the Feichtinger algebra increased with the turn of the century, primarily due to \cite{grochenig2013foundations, reiter2000}. Since then, the Feichtinger algebra and the more general modulation spaces have appeared in many books, e.g.\ \cite{modulation_spaces2020, time_frequency2020, pdebook2011}. We refer the reader to \cite{hans_look_back_paper} for interesting historical remarks and references to various applications of the Feichtinger algebra. For more recent applications we recommend \cite{Are2020, Bayer2020, berge2019modulation, Deformation15, Franz19, Sugimoto2015}.

\subsubsection*{Outline}
We begin in Section \ref{sec: 2} by motivating the need for the Feichtinger algebra $\mathbf{S}_{0}(\mathbb{R})$. In Section \ref{sec: 3} we introduce the Feichtinger algebra $\mathbf{S}_{0}(\mathbb{R})$ on the real line. Moreover, we connect the Feichtinger algebra with the Heisenberg group and the short-time Fourier transform. In Section \ref{sec: 4} we discuss the various properties of $\mathbf{S}_{0}(\mathbb{R})$. We give an alternative definition of $\mathbf{S}_{0}(\mathbb{R})$ in Section \ref{sec: 5} based on a frequency decomposition. Finally, in Section \ref{sec: 6} we outline some generalizations to give the reader a more holistic view. 

\subsubsection*{Acknowledgements}
I would like to thank Stine Marie Berge and Franz Luef for helpful suggestions and corrections. This note is not intended for publication as it does not contain an ounce of originality from my part. On the positive side, this allows the note to be rapidly updated in case of errors or wrongful omissions. If any such issues are encountered when reading, please contact me at \texttt{eirik.berge@ntnu.no} so that these can get fixed.

\section{Where Lebesgue and Schwartz Come Up Short}
\label{sec: 2}

To properly motivate the Feichtinger algebra $\mathbf{S}_{0}(\mathbb{R})$ on the real line, we will outline some problems with the traditional spaces $L^{2}(\mathbb{R})$ and $\mathscr{S}(\mathbb{R})$. The notation $L^{p}(\mathbb{R})$ for $1 \leq p \leq \infty$ indicates the space of measurable functions $f:\mathbb{R} \to \mathbb{C}$ satisfying
\begin{equation}
\label{eq:lebesgue_norm}
    \|f\|_{L^{p}(\mathbb{R})} \coloneqq \left(\int_{-\infty}^{\infty}|f(t)|^{p} \, dt\right)^{\frac{1}{p}} < \infty.
\end{equation}
For $p = \infty$ we exchange the expression \eqref{eq:lebesgue_norm} by the essential supremum. We conceal the fact that elements in $L^{p}(\mathbb{R})$ are really equivalence classes to keep our sanity intact. The notation $\mathscr{S}(\mathbb{R})$ denotes the \textit{Schwartz functions}, that is, the smooth and rapidly decaying functions on the real line. Both spaces $L^{2}(\mathbb{R})$ and $\mathscr{S}(\mathbb{R})$ are invaluable in many areas of analysis, e.g.\ Fourier analysis and PDE theory. \par
Time-frequency analysis is often concerned with a function $f$ and its frequency information given by $\mathcal{F}f$, where $\mathcal{F}$ denotes the \textit{Fourier transform} defined by 
\begin{equation}
\label{eq:fourier_transform}
    \mathcal{F}f(\omega) \coloneqq \int_{-\infty}^{\infty}f(t)e^{-2\pi i \omega t} \, dt.
\end{equation}
In the case where $f$ represents an audio signal, then $f(t)$ is the amplitude at time $t \in \mathbb{R}$ and $\mathcal{F}f(\omega)$ represents the amplitude of the frequency $\omega \in \mathbb{R}$ in the signal $f$. So why are the spaces $L^{2}(\mathbb{R})$ and $\mathscr{S}(\mathbb{R})$ used in time-frequency analysis? First of all, both $L^{2}(\mathbb{R})$ and $\mathscr{S}(\mathbb{R})$ are preserved under the Fourier transform. For $L^{2}(\mathbb{R})$, the Fourier transform $\mathcal{F}$ is defined by \eqref{eq:fourier_transform} on the dense subspace $\mathscr{S}(\mathbb{R})$ and then continuously extended to $L^{2}(\mathbb{R})$. Secondly, both $L^{2}(\mathbb{R})$ and $\mathscr{S}(\mathbb{R})$ are invariant under \textit{time-shifts} $T_{x}$ for $x \in \mathbb{R}$ and \textit{frequency-shifts} $M_{\omega}$ for $\omega \in \mathbb{R}$, where
\[T_{x}f(t) \coloneqq f(t - x), \qquad M_{\omega}f(t) \coloneqq e^{2 \pi i \omega t}f(t).\] 
Although this is a promising start, the spaces $L^{2}(\mathbb{R})$ and $\mathscr{S}(\mathbb{R})$ have problems that we now discuss.

\subsubsection*{Square Integrable Functions}

The space $L^{2}(\mathbb{R})$ is well-behaved as a collection since it is a Hilbert space. Nevertheless, individual elements can be very pathological; a generic function in $L^{2}(\mathbb{R})$ does not even have a continuous representative. 

\begin{example}
One would naively expect that if $f \in L^{2}(\mathbb{R})$ is sufficiently smooth, then the samples $(f(n))_{n \in \mathbb{Z}}$ decay nicely. However, consider for $s \geq 0$ the awful function
\[f^{s}(t) \coloneqq \sum_{n = -\infty}^{\infty}\eta_{A_n}^{s}(t), \qquad A_{n} \coloneqq \left[n -\frac{1}{2^{|n| + 1}}, n + \frac{1}{2^{|n| + 1}}\right],\]
where $\eta_{A_{n}}^{s}$ is a smooth function satisfying $\textrm{supp}(\eta_{A_{n}}^{s}) \subset A_{n}$ and \[\max_{t \in A_{n}}|\eta_{A_{n}}^{s}(t)| = |\eta_{A_{n}}^{s}(n)| = |n|^{s}.\] 
One now checks, with mild revulsion, that $f^{s} \in L^{2}(\mathbb{R})$ is a smooth function satisfying \[f^{s}(n) = |n|^{s}, \qquad n \in \mathbb{Z}.\]
\end{example}

In summary, although $L^{2}(\mathbb{R})$ is collectively well-behaved, the individual elements in $L^{2}(\mathbb{R})$ can be very pathological.
    
\subsubsection*{Schwartz Functions}

The Schwartz functions $\mathscr{S}(\mathbb{R})$ are individually so well-behaved that the pointwise product of $f,g \in \mathscr{S}(\mathbb{R})$ satisfies $f \cdot g \in \mathscr{S}(\mathbb{R})$. However, the individual niceness comes at the cost of weak collective properties. There is unfortunately no natural way to make $\mathscr{S}(\mathbb{R})$ into a Banach space. Hence one has to make do with the Fréchet structure on $\mathscr{S}(\mathbb{R})$ induced by a countable collection of semi-norms. Be wary of anyone claiming that countably many semi-norms are \textquote{just as good as a single norm}. They are not. There is no way to make the \textit{tempered distributions} $\mathscr{S}'(\mathbb{R})$, the space of anti-linear continuous functionals on $\mathscr{S}(\mathbb{R})$, into a Fréchet space. Additionally, many functions of practical interest such as 
\begin{equation}
\label{eq:not_rapid_decay}
    f_{a}(t) \coloneqq \frac{1}{a^2 + t^2}, \quad a > 0,
\end{equation}
are simply not nice enough to be in $\mathscr{S}(\mathbb{R})$. In summary, although individual elements in $\mathscr{S}(\mathbb{R})$ are well-behaved, the space  $\mathscr{S}(\mathbb{R})$ has poor collective properties.

\section{Introducing the Feichtinger Algebra}
\label{sec: 3}

In Section \ref{sec: 2} we saw that there is a trade-off between the collectively well-behaved $L^{2}(\mathbb{R})$ and the individually well-behaved $\mathscr{S}(\mathbb{R})$. The Feichtinger algebra provides a sweet spot between these two extremes. 

\begin{definition}
We define the \textit{Feichtinger algebra} $\mathbf{S}_{0}(\mathbb{R})$ as all elements $f \in L^{2}(\mathbb{R})$ satisfying
\begin{equation}
    \label{eq:explicit_Feichtinger_norm}
    \|f\|_{\mathbf{S}_{0}(\mathbb{R})} \coloneqq \int_{-\infty}^{\infty}\int_{-\infty}^{\infty} \left|\int_{-\infty}^{\infty}f(t)e^{-\pi(x-t)^2}e^{-2\pi i \omega t}\, dt\right| \, dx \, d\omega < \infty.
\end{equation}
\end{definition}

It is hard to muster any initial enthusiasm when looking at \eqref{eq:explicit_Feichtinger_norm}. The expression somehow manages to be both intimidating and seemingly random. It is not even obvious that \eqref{eq:explicit_Feichtinger_norm} actually defines a norm on $\mathbf{S}_{0}(\mathbb{R})$. We now demystify the Feichtinger algebra $\mathbf{S}_{0}(\mathbb{R})$ by linking it to the representation theory of the Heisenberg group.

\subsubsection*{Motivation Through the Heisenberg Group}

The three-dimensional \textit{Heisenberg group} $\mathbb{H}$ is an important structure in topics ranging from sub-Riemannian geometry \cite[Section 4.4.3]{subRiemannian2019} to quantum mechanics \cite[Chapter 14]{hall2013quantum}. The ubiquity of the Heisenberg group in modern mathematics is demonstrated in \cite{howe1980role}. As a set, the Heisenberg group is simply $\mathbb{R}^{3}$ with the group multiplication
\[\big(x,\omega,\tau \big) \cdot \big(x',\omega',\tau'\big) := \left(x + x', \omega + \omega', \tau + \tau' + \frac{1}{2}(x' \omega - x \omega')\right).\]
The popularity of the Heisenberg group can partly be explained by the fact that its representation theory is completely understood. The famous Stone-von Neumann theorem, see e.g.\ \cite[Section 9.3]{grchenig2020new}, shows that the most important representation of the Heisenberg group $\mathbb{H}$ is the \textit{Schr\"{o}dinger representation} $\rho$ acting on $f \in L^{2}(\mathbb{R})$ by 
\[\rho(x, \omega, \tau)f(t) \coloneqq e^{2 \pi i \tau}e^{\pi i x \omega}T_{x}M_{\omega}f(t), \qquad (x, \omega, \tau) \in \mathbb{H}.\]
For those familiar with representation theory, the Schr\"{o}dinger representation $\rho$ is an irreducible unitary representation. \par
Let us define the \textit{short-time Fourier transform} $V_{g}f$ of $f \in L^{2}(\mathbb{R})$ with respect to $g \in L^{2}(\mathbb{R})$ to be 
\begin{equation*}
V_{g}f(x, \omega) \coloneqq  \langle f, M_{\omega}T_{x}g \rangle = \int_{-\infty}^{\infty}f(t)\overline{g(t - x)}e^{-2 \pi i \omega t } \, dt.
\end{equation*}
The short-time Fourier transform is related to the Schr\"{o}dinger representation $\rho$ by the formula
\[\langle f, \rho(x, \omega, \tau)g \rangle = e^{-2 \pi i \tau} e^{\pi i x \omega} V_{g}f(x, \omega),\]
for $(x, \omega, \tau) \in \mathbb{H}$ and $f, g \in L^{2}(\mathbb{R})$. So what does all this have to do with the Feichtinger algebra $\mathbf{S}_{0}(\mathbb{R})$? By considering the Gaussian $g(t) \coloneqq e^{-\pi t^2}$ and comparing with \eqref{eq:explicit_Feichtinger_norm} we find that
\begin{equation}
\label{eq:STFT}
    \|f\|_{\mathbf{S}_{0}(\mathbb{R})} = \|V_{g}f\|_{L^{1}(\mathbb{R}^{2})}.
\end{equation}
Since $f \in \mathscr{S}(\mathbb{R})$ implies that $V_{g}f \in L^{1}(\mathbb{R}^2)$ by an elementary computation, it follows that \[\mathscr{S}(\mathbb{R}) \subset \mathbf{S}_{0}(\mathbb{R}) \subset L^{2}(\mathbb{R}).\] \par
The short-time Fourier transform $V_{g}f$ is, in the language of applied time-frequency analysis \cite[Chapter 3]{grochenig2013foundations}, a joint time-frequency depiction of $f$. 
As such, we should suspect that elements in $\mathbf{S}_{0}(\mathbb{R})$ are reasonably well-behaved both in time and frequency. Moreover, it should be noted that there is nothing special about choosing the Gaussian; any non-zero $g \in \mathscr{S}(\mathbb{R})$ can be used in \eqref{eq:STFT} to give an equivalent norm on $\mathbf{S}_{0}(\mathbb{R})$ by \cite[Proposition 4.10]{mads18}. Armed with this refreshing view of the Feichtinger algebra $\mathbf{S}_{0}(\mathbb{R})$, we can more easily develop its basic properties.

\section{Breaking the Trade-Off}
\label{sec: 4}

In this section we develop the basic properties of the Feichtinger algebra $\mathbf{S}_{0}(\mathbb{R})$. The crux is that $\mathbf{S}_{0}(\mathbb{R})$ is a well-behaved Banach space that lies in-between $\mathscr{S}(\mathbb{R})$ and $L^{2}(\mathbb{R})$.

\begin{description}
\item[Banach Space Structure:] It is well known \cite[Corollary 1.2.12]{time_frequency2020} that the short-time Fourier transform $V_{g}f$ is continuous and satisfies the norm equality
\begin{equation}
\label{eq:stft_orthogonality}
    \|V_{g}f\|_{L^{2}(\mathbb{R}^{2})} = \|f\|_{L^{2}(\mathbb{R})}\|g\|_{L^{2}(\mathbb{R})}.
\end{equation}
Hence if $\|V_{g}f\|_{L^{1}(\mathbb{R}^{2})} = 0$, then it follows from \eqref{eq:stft_orthogonality} that $f = 0$. This implies that \eqref{eq:explicit_Feichtinger_norm} is indeed a norm on $\mathbf{S}_{0}(\mathbb{R})$. The same conclusion can also be deduced from the irreducibility of the Schr\"{o}dinger representation. Showing completeness of $\mathbf{S}_{0}(\mathbb{R})$ with the norm $\|\cdot\|_{\mathbf{S}_{0}(\mathbb{R})}$ is relatively straightforward, see \cite[Theorem 4.12]{mads18} for details. 
\item[Time-Frequency Shift Invariance:] The short-time Fourier transform is easily seen to satisfy the properties \begin{align*}
    V_{g}(T_{y}f)(x, \omega) & = V_{g}f(x - y, \omega), \\ 
    V_{g}(M_{\eta}f)(x, \omega) & = V_{g}f(x, \omega - \eta),
\end{align*}
for $f,g \in L^{2}(\mathbb{R})$ and $x,y,\omega, \eta \in \mathbb{R}$. As such, it follows from \eqref{eq:STFT} that $\mathbf{S}_{0}(\mathbb{R})$ is invariant under time-shifts and frequency-shifts. The Feichtinger algebra $\mathbf{S}_{0}(\mathbb{R})$ is, in a sense made precise in \cite[Theorem 12.1.9]{grochenig2013foundations}, the smallest space suitably invariant under time-shifts and frequency-shifts.
\item[Fourier Invariance and Continuity:] The well-known Plancherel theorem states that the Fourier transform is a bijective isometry $\mathcal{F}:L^{2}(\mathbb{R}) \to L^{2}(\mathbb{R})$. This immediately implies the \textit{fundamental identity of time-frequency analysis}
\begin{equation}
\label{eq:fundamental_identity}
    V_{g}f(x, \omega) = e^{-2 \pi i x \omega}V_{\mathcal{F}g}\mathcal{F}f(\omega, -x).
\end{equation}
The fact that $\mathcal{F}g = g$ for the Gaussian $g(t) \coloneqq e^{-\pi t^{2}}$ shows together with \eqref{eq:fundamental_identity} that $\mathcal{F}:\mathbf{S}_{0}(\mathbb{R}) \to \mathbf{S}_{0}(\mathbb{R})$ is a bijective isometry as well. Since $\mathbf{S}_{0}(\mathbb{R})\subset L^{1}(\mathbb{R})$ follows from \cite[Corollary 4.2 (iii)]{mads18}, we conclude that both $f$ and $\mathcal{F}f$ are continuous when $f \in \mathbf{S}_{0}(\mathbb{R})$.
\end{description}

\begin{example}
\label{example:not_in_Feichtinger}
    The function \[f(t) \coloneqq \frac{e^{2 \pi i t} - 1}{2 \pi i t}\]
    is both continuous and in $L^{2}(\mathbb{R})$. In spite of this, it is not in $\mathbf{S}_{0}(\mathbb{R})$ since 
    \[\mathcal{F}f(\omega) = \chi_{[0, 1]}(\omega),\]
    which is not continuous.
\end{example}

The following result is taken from \cite[Proposition 12.1.6]{grochenig2013foundations} and shows that sufficient time-frequency decay implies membership in $\mathbf{S}_{0}(\mathbb{R})$.

\begin{lemma}
\label{lemma_decay}
    Assume that $f \in L^{2}(\mathbb{R})$ satisfies \[\left(\int_{-\infty}^{\infty}|f(t)|^{2}(1 + |t|)^{1 + \epsilon} \, dt\right)\left(\int_{-\infty}^{\infty}|\mathcal{F}f(\omega)|^{2}(1 + |\omega|)^{1 + \epsilon} \, d\omega\right) < \infty,\]
    for some $\epsilon > 0$. Then $f \in \mathbf{S}_{0}(\mathbb{R})$.
\end{lemma}
    
\begin{example}
\label{ex: non-rapid example}
Consider for $a > 0$ the function $f_{a}$ given in \eqref{eq:not_rapid_decay}. Verifying that $f_{a} \in \mathbf{S}_{0}(\mathbb{R})$ directly by using \eqref{eq:explicit_Feichtinger_norm} is a waste of a fine afternoon. However, a standard residue calculation shows that the Fourier transform of $f_{a}$ is given by
    \[\mathcal{F}f_{a}(\omega) =  \frac{\pi}{a}e^{-2\pi a |\omega|}.\]
We can now apply Lemma \ref{lemma_decay} to see that $f_{a} \in \mathbf{S}_{0}(\mathbb{R})$.
\end{example}
    
Even though $\mathscr{S}(\mathbb{R}) \neq \mathbf{S}_{0}(\mathbb{R})$ as demonstrated in Example \ref{ex: non-rapid example}, the Schwartz space $\mathscr{S}(\mathbb{R})$ is dense in $\mathbf{S}_{0}(\mathbb{R})$ by \cite[Proposition 11.3.4]{grochenig2013foundations}. The following result, see e.g.\ \cite[Corollary 4.14]{mads18}, explains why $\mathbf{S}_{0}(\mathbb{R})$ is called the Feichtinger algebra, and not the Feichtinger space.

\begin{proposition}
The pointwise product of $f, g \in \mathbf{S}_{0}(\mathbb{R})$ satisfies $f \cdot g \in \mathbf{S}_{0}(\mathbb{R})$. Hence the Fourier invariance of $\mathbf{S}_{0}(\mathbb{R})$ implies that the convolution
\[(f * g)(t) \coloneqq \int_{-\infty}^{\infty}f(s)g(t - s) \, ds\]
of $f,g \in \mathbf{S}_{0}(\mathbb{R})$ satisfies $f * g \in \mathbf{S}_{0}(\mathbb{R})$ as well.
\end{proposition}

Finally, let us consider the problem of sampling. We showed in Section \ref{sec: 2} that there are smooth functions $f \in L^{2}(\mathbb{R})$ where the samples $(f(n))_{n \in \mathbb{Z}}$ do not have any decay. For elements in $\mathbf{S}_{0}(\mathbb{R})$ this degeneracy is no longer possible as \cite[Theorem 5.7 (ii)]{mads18} shows. 

\begin{theorem}
\label{thm:restriction}
The samples $(f(n))_{n \in \mathbb{Z}}$ of $f \in \mathbf{S}_{0}(\mathbb{R})$ satisfy
\begin{equation}
\label{eq:restriction_equation}
    \sum_{n = -\infty}^{\infty}|f(n)| \leq C \cdot \|f\|_{\mathbf{S}_{0}(\mathbb{R})},
\end{equation}
for some absolute constant $C > 0$.
\end{theorem}

\begin{remark}
It follows from Example \ref{example:not_in_Feichtinger} that there are elements $f \in L^{2}(\mathbb{R}) \setminus \mathbf{S}_{0}(\mathbb{R})$ such that $f$ and $\mathcal{F}f$
satisfy \eqref{eq:restriction_equation}.
\end{remark}

There are other intriguing properties of the Feichtinger algebra $\mathbf{S}_{0}(\mathbb{R})$ that we do not go into. In particular, we encourage the reader to seek out the Poisson summation formula for $\mathbf{S}_{0}(\mathbb{R})$ \cite[Corollary 12.1.5]{grochenig2013foundations}, the tensor factorization property for $\mathbf{S}_{0}(\mathbb{R})$ \cite[Theorem 7.4]{mads18}, and the kernel theorem for $\mathbf{S}_{0}(\mathbb{R})$ \cite[Theorem 9.3]{mads18}. Moreover, if we let $\mathbf{S}_{0}'(\mathbb{R})$ denote the dual space of $\mathbf{S}_{0}(\mathbb{R})$, then $\mathbf{S}_{0}(\mathbb{R}) \subset L^{2}(\mathbb{R}) \subset \mathbf{S}_{0}'(\mathbb{R})$ is a \textit{Banach Gelfand triple} \cite{Feichtinger2008}.

\section{A Frequency Decomposition}
\label{sec: 5}
We briefly outline an alternative definition of the Feichtinger algebra $\mathbf{S}_{0}(\mathbb{R})$ based on a frequency decomposition. This viewpoint is more geometric and is often preferred in the PDE literature, see e.g.\ \cite{pdebook2011}.
We first decompose $\mathbb{R}$, viewed as the frequency domain, into discrete containers. Consider the covering $(Q_{n})_{n \in \mathbb{Z}}$ given by \[\mathbb{R} = \bigcup_{n = -\infty}^{\infty} Q_{n} \coloneqq \bigcup_{n = -\infty}^{\infty} (n, n + 2).\] The goal is to associate to $f \in L^{2}(\mathbb{R})$ a collection of numbers $A_{n}(f)$ for $n \in \mathbb{Z}$. The number $A_{n}(f)$ should represent the average amplitude of the frequencies of $f$ in the container $Q_{n}$. 

\begin{remark}
It is tempting to define the numbers $(A_{n}(f))_{n \in \mathbb{Z}}$ for $f \in L^{2}(\mathbb{R})$ as \begin{equation}
\label{eq:wrong_frequency_def}
    A_{n}(f) \coloneqq \frac{1}{2}\|\mathcal{F}^{-1}(\chi_{Q_{n}} \cdot \mathcal{F}f)\|_{L^{1}(\mathbb{R})} =  \frac{1}{2}\int_{-\infty}^{\infty}\left|\left(\mathcal{F}^{-1}(\chi_{Q_{n}} \cdot \mathcal{F}f)\right)(t)\right| \, dt.
\end{equation}
However, consider the case where $f \in L^{2}(\mathbb{R})$ satisfies $\mathcal{F}f \equiv 1$ on $Q_{n}$ for some $n \in \mathbb{Z}$. Then \eqref{eq:wrong_frequency_def} fails spectacularly at capturing the average amplitude of the frequencies of $f$ in the container $Q_{n}$ since  
\[A_{n}(f) = \frac{1}{2}\int_{-\infty}^{\infty}\left|\left(\mathcal{F}^{-1}\chi_{Q_{n}}\right)(t)\right| \, dt = \frac{1}{4\pi}\int_{-\infty}^{\infty}\left|\frac{e^{4 \pi i t} - 1}{t}\right| \, dt = \infty.\]
The problem lies with the abrupt cutoffs in the characteristic functions $(\chi_{Q_{n}})_{n \in \mathbb{Z}}$.
\end{remark}

To remedy the problem in the remark above, we consider a collection $(\psi_{n})_{n \in \mathbb{Z}}$ of smooth functions $\psi_{n}:\mathbb{R} \to [0, 1]$ satisfying \[\textrm{supp}(\psi_{n}) \subset Q_{n}, \qquad \sum_{n= -\infty}^{\infty}\psi_{n}(t) = 1, \qquad \sup_{n \in \mathbb{Z}}\int_{-\infty}^{\infty}\left|\left(\mathcal{F}^{-1}\psi_{n}\right)(t)\right| \, dt < \infty.\] The existence of $(\psi_{n})_{n \in \mathbb{Z}}$ is guaranteed by e.g.\ \cite[Chapter 6.2]{pdebook2011} and is only possible due to the overlap $Q_{n} \cap Q_{k} \neq \emptyset$ whenever $k \in \{n-1, n, n+1\}$. We can now for $f \in L^{2}(\mathbb{R})$ properly define the numbers $(A_{n}(f))_{n \in \mathbb{Z}}$ as
\begin{equation*}
    A_{n}(f) \coloneqq \frac{1}{2}\|\mathcal{F}^{-1}(\psi_{n} \cdot \mathcal{F}f)\|_{L^{1}(\mathbb{R})} = \frac{1}{2}\int_{-\infty}^{\infty}\left|\left(\mathcal{F}^{-1}(\psi_{n} \cdot \mathcal{F}f)\right)(t)\right| \, dt.
\end{equation*}
The following result originates from \cite{modulation_space_report} and gives an alternative definition of $\mathbf{S}_{0}(\mathbb{R})$.

\begin{theorem}
The Feichtinger algebra $\mathbf{S}_{0}(\mathbb{R})$ consists precisely of the elements $f \in L^{2}(\mathbb{R})$ where the numbers $(A_{n}(f))_{n \in \mathbb{Z}}$ satisfy
\begin{equation}
\label{eq:freq_decomposition_def_Feichtinger}
    \sum_{n = -\infty}^{\infty}A_{n}(f) < \infty.
\end{equation}
In fact, the expressions \eqref{eq:explicit_Feichtinger_norm} and \eqref{eq:freq_decomposition_def_Feichtinger} are equivalent norms on $\mathbf{S}_{0}(\mathbb{R})$.
\end{theorem}

\section{Various Generalizations}
\label{sec: 6}

Now that we have demonstrated the impressive features of the Feichtinger algebra $\mathbf{S}_{0}(\mathbb{R})$, it is time to point out different generalizations. By using the connection with the short-time Fourier transform in \eqref{eq:STFT}, we have the following obvious generalization.

\begin{definition}
For $1 \leq p \leq \infty$ the space $M^{p}(\mathbb{R})$ consists of the tempered distributions $f \in \mathscr{S}'(\mathbb{R})$ such that \begin{equation}
\label{eq:modulation_definition}
    \|f\|_{M^{p}(\mathbb{R})} \coloneqq \|V_{g}f\|_{L^{p}(\mathbb{R}^{2})} < \infty,
\end{equation}
where $g(t) \coloneqq e^{-\pi t^2}$ is the Gaussian.
\end{definition}

\begin{remark}
The reader might be weary that we are applying the short-time Fourier transform to elements in $\mathscr{S}'(\mathbb{R})$ that are not in $L^{2}(\mathbb{R})$. However, this is perfectly valid by interpreting $V_{g}f(x, \omega) = \langle f, M_{\omega}T_{x}g \rangle$ as the duality pairing between $\mathscr{S}'(\mathbb{R})$ and $\mathscr{S}(\mathbb{R})$. 
\end{remark}

In our new notation we have $M^{1}(\mathbb{R}) = \mathbf{S}_{0}(\mathbb{R})$ and $M^{\infty}(\mathbb{R}) = \mathbf{S}_{0}'(\mathbb{R})$. As with the Feichtinger algebra, the choice to consider the Gaussian is out of convenience rather than necessity. It should be noted that $M^{p}(\mathbb{R}) \subset L^{2}(\mathbb{R})$ only for $1 \leq p  \leq 2$ with $M^{2}(\mathbb{R}) = L^{2}(\mathbb{R})$. Moreover, one has for $1\leq p \leq q \leq \infty$ the general inclusion relation \[M^{p}(\mathbb{R}) \subset M^{q}(\mathbb{R}).\] A straightforward computation shows that $\chi_{[0, 1]} \in M^{p}(\mathbb{R})$ for all $p > 1$. Hence the fact that elements in $M^{1}(\mathbb{R})$ are continuous does not generalize to any of the other modulation spaces. More basic properties of modulation spaces can be found in e.g.\ \cite[Chapter 2]{time_frequency2020}. 

\begin{example}
The famous \textit{Dirac Comb distribution} $\delta_{\mathbb{Z}} \in \mathscr{S}'(\mathbb{R})$ is defined as acting on $f \in \mathscr{S}(\mathbb{R})$ by
\begin{equation*}
    \delta_{\mathbb{Z}}(f) := \sum_{n = -\infty}^{\infty}f(n).
\end{equation*}
Theorem \ref{thm:restriction} gives that $\delta_{\mathbb{Z}} \in M^{\infty}(\mathbb{R})$. 
\end{example}

There are three minor generalizations of the modulation spaces $M^{p}(\mathbb{R})$ that can immediately be considered:

\begin{description}
\item[Increasing the dimension:] We have restricted ourselves to modulation spaces $M^{p}(\mathbb{R})$ on $\mathbb{R}$, but this is only for simplicity. One can consider the \textit{$n$-dimensional modulation spaces} $M^{p}(\mathbb{R}^{n})$ for $1 \leq p \leq \infty$ with only minor modifications.
\item[Mixed Integrability:] In \eqref{eq:modulation_definition} we measured the size of $V_{g}f$ with the $L^{p}$-norm. One can consider the \textit{mixed spaces} $L^{p,q}$ for $1 \leq p,q \leq \infty$ where 
\begin{equation}
\label{eq:mixed_modulation_spaces}
    \|V_{g}f\|_{L^{p,q}(\mathbb{R}^{2})} \coloneqq \left(\int_{-\infty}^{\infty}\left(\int_{-\infty}^{\infty}|V_{g}f(x, \omega)|^{p}\, dx\right)^{\frac{q}{p}} \, d\omega \right)^{\frac{1}{q}}.
\end{equation}
When $p = \infty$ or $q = \infty$ we use the appropriate supremums in \eqref{eq:mixed_modulation_spaces}. This gives rise to the \textit{mixed modulation spaces} $M^{p,q}(\mathbb{R})$ for $1 \leq p,q \leq \infty$. Notice that when $p = q$ we regain the space $M^{p}(\mathbb{R}) = M^{p,p}(\mathbb{R})$.
\item[Introducing Weights: ] One can introduce weights into the mix to get more general spaces. For concreteness, let us consider the family of weights \[v_{s}(x,\omega) := (1 + |x|^2 + |\omega|^2)^{\frac{s}{2}}, \quad (x,\omega) \in \mathbb{R}^{2}, \, s \geq 0.\] The family $v_s$ is sometimes referred to as the \textit{polynomial weights}. Using these weights, one obtains the \textit{polynomial modulation spaces} $M_{s}^{p}(\mathbb{R})$ for $1 \leq p \leq \infty$ and $s \geq 0$ where 
\[\|f\|_{M_{s}^{p}(\mathbb{R})} \coloneqq \|V_{g}f \cdot v_{s}\|_{L^{p}(\mathbb{R}^{2})}.\]
It follows from \cite[Proposition 11.3.1]{grochenig2013foundations} that the \textit{Bessel potential spaces} \[H^{s}(\mathbb{R}) = \left\{f \in \mathscr{S}'(\mathbb{R}) \, : \, \int_{-\infty}^{\infty}|\mathcal{F}f(\omega)|^2 (1 + |\omega|^2)^s \, d\omega < \infty\right\}, \qquad s \geq 0,\] 
can be represented as weighted modulation spaces.
\end{description}

In practice, all three generalizations above can be considered at once. Hence one often see authors juggling the \textit{general modulation spaces} $M_{w}^{p,q}(\mathbb{R}^{n})$, where $1 \leq p,q \leq \infty$ and $w$ a suitable weight function. Although this gives a general and powerful family of spaces, it should not be forgotten that $p = q = 1$ and $w \equiv 1$ gives a particularly nice space; the \textit{$n$-dimensional Feichtinger algebra} $\mathbf{S}_{0}(\mathbb{R}^{n})$. The results in Section \ref{sec: 4} and Section \ref{sec: 5} are still valid for the $n$-dimensional Feichtinger algebra $\mathbf{S}_{0}(\mathbb{R}^{n})$ with the obvious modifications. \par 
The general modulation spaces $M_{w}^{p,q}(\mathbb{R}^{n})$ are concrete examples of \textit{coorbit spaces}. Said briefly, coorbit spaces are function spaces that are constructed out of well-behaved representations of locally compact groups. We refer the reader to the original paper \cite{feichtinger1988unified} and the recent survey \cite{berge2021coorbitsurvey} for more on coorbit spaces. 

\subsubsection*{Two Further Generalizations}
For readers with curious minds and strong stomachs, there are two broader generalizations of the Feichtinger algebra we mention:
\begin{itemize}
    \item The class of locally compact abelian groups includes the Euclidean spaces $\mathbb{R}^{n}$, the integer lattices $\mathbb{Z}^{n}$, and the tori $\mathbb{T}^{n}$. The Feichtinger algebra $\mathbf{S}_{0}(G)$ was originally defined in \cite{feichtinger_original} for a locally compact abelian group $G$. This is also the setting where the survey \cite{mads18} takes place. Most properties we have discussed carries over in some way to locally compact abelian groups. However, working in this setting introduces significant technicalities and many genuine difficulties.
    \item Recently, the Feichtinger algebra has been extended to certain nilpotent Lie groups in \cite{berge2019modulation} based on the tour de force paper \cite{fischer2018heisenberg}. Rather than using a generalization of the short-time Fourier transform, the approach uses frequency decompositions motivated by Section~\ref{sec: 5}. We refer the reader to \cite{Hans_2, Hans_Grobner} for the origins of the general \textit{decomposition spaces}. Extending time-frequency analysis to the nilpotent setting has gained traction in later years, see e.g.\  \cite{grchenig2020new} for a recent interesting example.
\end{itemize}

For those interested in learning more about the Feichtinger algebra, the well-written book \cite[Chapter 12]{grochenig2013foundations} and the enjoyable survey \cite{mads18} are recommended as good starting points.

\bibliographystyle{abbrv}
\bibliography{main}

\end{document}